\documentclass[11pt]{article}%
\usepackage{a4}%
\usepackage[centertags]{amsmath}%
\usepackage{eucal,dsfont,accents,bbm,amsfonts,amssymb,amsthm,amsopn}%
\usepackage[usenames]{color}
\definecolor{citegreen}{rgb}{0,0.6,0}
\definecolor{refred}{rgb}{0.8,0,0}
\usepackage[colorlinks, citecolor=citegreen, linkcolor=refred]{hyperref}
\title{A note on the compactness theorem for 4d Ricci shrinkers}
\author{Robert Haslhofer, Reto M\"{u}ller}
\date{}
\providecommand{\abs}[1]{\lvert #1\rvert}%
\DeclareMathOperator{\Hess}{Hess}%
\DeclareMathOperator{\Vol}{Vol}%
\DeclareMathOperator{\Rm}{Rm}%
\DeclareMathOperator{\Rc}{Rc}%

\newcommand{\RR}{\mathbb{R}}%
\newcommand{\Lap}{\triangle}%
\newcommand{\D}{\nabla}%
\newcommand{\sW}{\mathcal{W}}%
\newcommand{\ol}{\overline}%
\newcommand{\ul}{\underline}%
\newtheoremstyle{break}%
  {12pt}%
  {16pt}%
  {\itshape}%
  {}%
  {\bfseries}%
  {}%
  {\newline}%
  {\thmname{#1}\thmnumber{ #2}\thmnote{ \normalfont{(#3)}}}%

\theoremstyle{definition}%
\theoremstyle{remark}%
\newtheorem*{rem}{Remark}%
\theoremstyle{break}%
\newtheorem{lemma}{Lemma}[section]%
\newtheorem{thm}[lemma]{Theorem}%
\numberwithin{equation}{section}%
\begin{document}%
\maketitle%
\pagenumbering{arabic}%
\begin{abstract}
In \cite{HM11} we proved an orbifold Cheeger-Gromov compactness theorem for complete 4d Ricci shrinkers
with a lower bound for the entropy, an upper bound for the Euler characterisic, and a lower bound for the gradient of the potential at large distances.
In this note, we show that the last two assumptions in fact can be removed.
The key ingredient is a recent estimate of Cheeger-Naber \cite{CN14}.
\end{abstract}

\section{Introduction}
The goal of this short note is to improve our compactness theorem for 4d Ricci shrinkers from \cite{HM11}, by removing two of the three assumptions.
Recall that a \emph{Ricci shrinker} is a Riemannian manifold $(M,g)$ (smooth, complete, connected) together with a smooth function $f:M\to\RR$ such that
\begin{equation}\label{1.solitoneq}
\Rc_{g} + \Hess_{g} f = \tfrac{1}{2}g.
\end{equation}
Ricci shrinkers arise as singularity models for Hamilton's Ricci flow \cite{Ham95s}, and have received a lot of attention, especially in the last 10 years; see \cite{Cao09} for a recent survey.

As explained in \cite{HM11}, we can always find a natural basepoint $p\in M$ where the
potential $f$ attains its minimum, and we can always normalize $f$ such that
\begin{equation}\label{1.normalizationf}
\int_M (4\pi)^{-n/2}e^{-f}dV_g = 1.
\end{equation}
After imposing \eqref{1.normalizationf} every Ricci shrinker has a well defined \emph{Perelman entropy} \cite{Per02},
\begin{equation}\label{1.defmu}
\mu(g)= \sW(g,f) = \int_M\big(\abs{\D f}_g^2 + R_g + f -
n\big)(4\pi)^{-n/2}e^{-f}dV_g>-\infty.
\end{equation}

Let us now recall our compactness theorem for 4d Ricci shrinkers \cite[Thm 1.2]{HM11}.
We proved that any sequence $(M_i,g_i,f_i,p_i)$ of 4d Ricci shrinkers (with
normalization and basepoint $p_i$ as above) with entropy uniformly bounded below,
\begin{equation}\label{ass1}
\mu(g_i)\geq\ul{\mu}>-\infty,
\end{equation}
Euler-characterisic uniformly bounded above,
\begin{equation}\label{ass2}
\chi(M_i)\leq\bar{\chi}<\infty,
\end{equation}
and gradient of the potential uniformly bounded below at large distances,
\begin{equation}\label{ass3}
\abs{\D f_i}(x)\geq c>0 \qquad\text{if } d(x,p_i)\geq c^{-1},
\end{equation}
has a subsequence that converges to an orbifold Ricci shrinker in the pointed orbifold Cheeger-Gromov sense.
This means in particular that the convergence is smooth away from the orbifold singularities, which are isolated singularities modelled on $\RR^n/\Gamma$ for some finite subgroup $\Gamma\subset \mathrm{O}(n)$, see \cite[Sec. 3]{HM11} for the precise definitions.\\

The purpose of this short note is to improve our compactness theorem by removing the assumptions \eqref{ass2} and \eqref{ass3}. In other words, we prove the following theorem.

\begin{thm}[Compactness theorem for 4d Ricci shrinkers]\label{main_thm}
Let $(M_i,g_i,f_i)$ be a sequence of 4d Ricci shrinkers with entropy uniformly
bounded below, i.e. $\mu(g_i)\geq\ul{\mu}>-\infty$.
Then a subsequence of $(M_i,g_i,f_i,p_i)$ converges to an orbifold
Ricci shrinker in the pointed orbifold Cheeger-Gromov sense.
\end{thm}

Theorem \ref{main_thm} gives orbifold compactness of the space of (shrinking) singularity models for 4d Ricci flow, assuming only a lower bound for the entropy.

\begin{rem}
 Simple examples, like cylinders over $3$-dimensional lense spaces, show that the entropy assumption in Theorem \ref{main_thm} is indeed necessary.
 Due to Perelman's monotonicity formula \cite{Per02}, the entropy assumption is of course perfectly natural.
\end{rem}

Related interesting compactness theorems for Ricci solitons have been proved by Cao-Sesum \cite{CS07}, Weber \cite{Web08}, X.~Zhang \cite{Zha06}, Z.~Zhang \cite{Z10}, Tian-Zhang \cite{TZ12} and Chen-Wang \cite{CW12}.
One key feature that distinguishes Theorem \ref{main_thm} from these compactness theorems, in addition to removing some other assumptions, is that it applies in the setting of complete (possibly noncompact) manifolds. This is in fact crucial, since most interesting singularity models for the Ricci flow are noncompact.\\

The main ingredient in our proof is a recent estimate of Cheeger-Naber \cite{CN14}, that gives $L^2$-control for the Riemann tensor of noncollapsed 4d metrics with bounded Ricci curvature.
Though not directly applicable in our context, we can in fact combine this $L^2$-estimate with our previous compactness theorem \cite[Thm. 1.1]{HM11}
by making use of a nice observation of Z.~Zhang \cite{Z10} as well as some further uniform estimates from \cite{HM11}.

\section{The proof}\label{sec_proof}

Our previous proof of the 4d compactness theorem was based on a localized Gauss-Bonnet argument on 4d Ricci shrinkers \cite[Sec. 4]{HM11}, which -- under the assumptions \eqref{ass2} and \eqref{ass3} --
gave us the necessary local $L^2$ Riemann bounds to apply our general $n$-dimensional compactness theorem \cite[Thm 1.1]{HM11}.
We will now replace this localized Gauss-Bonnet argument by a new $L^2$-estimate, that works even without the assumptions \eqref{ass2} and \eqref{ass3}.\footnote{The estimate of Cheeger-Naber \cite[Thm. 1.5.]{CN14} -- which is the key ingredient -- is of course based on a local Gauss-Bonnet argument. However, the estimate is strong enough that it even gives topological control, which ultimately works as substitute for the assumptions \eqref{ass2} and \eqref{ass3}.}
Namely, we will prove that for any 4d Ricci shrinker $(M,g,f)$ (with normalization and basepoint $p$ as before) we have the $L^2$ curvature bound
\begin{equation}\label{eq_toshow}
 \int_{B_r(p)}\abs{\Rm}_g^2\; dV_g\leq C(r)
\end{equation}
for some universal function $C(r)=C_{\underline{\mu}}(r)<\infty$ depending only on a lower bound $\underline{\mu}$ for the entropy $\mu(g)$.\\

We start by recalling three lemmas from \cite{HM11}, see also Cao-Zhou \cite{CZ09} and Perelman \cite{Per02}.
First, by \cite[Lem. 2.1]{HM11} the potential $f$ satisfies the estimate
\begin{equation}\label{lemma1}
\tfrac{1}{4}\big(d(x,p)-20\big)_{\!+}^{2} \leq f(x)-\mu(g) \leq
\tfrac{1}{4}\big(d(x,p)+20\big)^2
\end{equation}
for all $x\in M$, where $a_+ := \max\{0,a\}$.
Second, by \cite[Lem. 2.2]{HM11} we have the volume growth estimate
\begin{equation}\label{2.ballseqn}
\Vol B_r(p)\leq \omega r^4,
\end{equation}
for some universal constant $\omega<\infty$.
Third, by \cite[Lem. 2.3]{HM11} there exists a function
$\kappa(r)=\kappa_{\ul{\mu}}(r)>0$ such that we have
the lower volume bound
\begin{equation}\label{noncoll}
\Vol B_\delta(x)\geq\kappa(r)\delta^4
\end{equation}
for every ball $B_\delta(x)\subset B_r(p)$, $0<\delta\leq 1$.\\

We also recall that the Bianchi identity forces the quantity $R+\abs{\D f}^2-f$ to be constant, in fact
\begin{equation}
 R+\abs{\D f}^2-f=-\mu(g),
\end{equation}
see \cite[(2.16)]{HM11}. In particular, by \eqref{lemma1} we have the upper bound
\begin{equation}\label{uppersum}
 R+\abs{\D f}^2\leq \tfrac{1}{4}\big(d(x,p)+20\big)^2.
\end{equation}
Since Ricci shrinkers always have nonnegative scalar curvature \cite{Zha09},
\begin{equation}\label{posscal}
R\geq 0,
\end{equation}
the estimate \eqref{uppersum} shows that both $R$ and $\abs{\D f}^2$ grow at most quadratically.
Finally, we recall that the entropy is automatically bounded from above,
\begin{equation}
\mu(g)\leq \ol{\mu},
\end{equation}
where $\ol{\mu}=\ol{\mu}_{\ul{\mu}}<\infty$ (and most likely $\ol{\mu}\leq 0$), see \cite[p. 1097]{HM11}.\\

Now, as in Z.~Zhang \cite{Z10} and Tian-Zhang \cite{TZ12} we consider the conformally rescaled metric
\begin{equation}
\tilde{g}=e^{-f}g.
\end{equation}
By the growth estimate \eqref{lemma1} the metrics are locally uniformly equivalent, i.e.
\begin{equation}\label{equivmet}
 \lambda^{-1}(r)g\leq \tilde{g}\leq \lambda(r)g \quad \textrm{on}\quad B_r(p),
\end{equation}
for some $\lambda(r)=\lambda_{\ul{\mu}}(r)<\infty$.
By the formula for the conformal transformation of the Ricci tensor, see e.g. \cite{Besse}, we have
\begin{align}
 \widetilde{\Rc}&=\Rc+\Hess f +\tfrac{1}{2}\nabla f\otimes \nabla f+\tfrac{1}{2}(\Lap f-\abs{\D f}^2)g\nonumber\\
&=\tfrac12 g+\tfrac{1}{2}\nabla f\otimes \nabla f+\tfrac{1}{2}(2-R-\abs{\D f}^2)g,
\end{align}
where we also used the soliton equation \eqref{1.solitoneq} and its trace.
Combining this with \eqref{uppersum}, \eqref{posscal} and \eqref{equivmet} we obtain the estimate
\begin{equation}\label{cn1}
 \abs{\widetilde{\Rc}}_{\tilde{g}}\leq K(r) \quad \textrm{on}\quad B_r(p),
\end{equation}
for some $K(r)=K_{\ul{\mu}}(r)<\infty$.\\

We will now cover $B_{r/2}(p)$ with suitable balls $B^{\tilde{g}}_{\bar{\delta}}(q_i)$ with center points $q_i\in B_{3r/4}(p)$.
Here, the notation $B^{\tilde{g}}$ indicates that the ball is defined with respect to the metric $\tilde{g}$.
We remark that one has to select the covering somewhat carefully, since in the noncompact case the manifold $(M,\tilde{g})$ is always incomplete.
However, if we chose $\bar{\delta}$ small enough to ensure that the balls $B^{\tilde{g}}_{2\bar{\delta}}(q_i)$ are contained in $B_r(p)$ then we are fine.
More precisely, using \eqref{2.ballseqn}, \eqref{noncoll} and \eqref{equivmet} we see that there exist constants $N(r)=N_{\ul{\mu}}(r)<\infty$, $v(r)=v_{\ul{\mu}}(r)>0$ and 
$\bar{\delta}(r)=\bar{\delta}_{\ul{\mu}}(r)\in (0,1]$ with the following properties. We can find $N(r)$ points $q_i\in B_{3r/4}(p)$ such that
\begin{align}
B^{\tilde{g}}_{2\bar{\delta}(r)}(q_i)\subseteq B_r(p),\label{cn2}\\
 \Vol^{\tilde{g}}(B^{\tilde{g}}_{\bar{\delta}(r)}(q_i))\geq v(r),\label{cn3}
\end{align}
and
\begin{equation}\label{eq_cov}
 B_{r/2}(p)\subseteq \bigcup_i B_{\bar{\delta}(r)}^{\tilde{g}}(q_i).
\end{equation}
By \eqref{cn1}, \eqref{cn2} and \eqref{cn3} we can now apply the fundamental estimate of Cheeger-Naber \cite[Thm. 1.5]{CN14}, which says that
\begin{equation}
 \int_{B^{\tilde{g}}_{\bar{\delta}(r)}(q_i)}\abs{\widetilde{\Rm}}_{\tilde{g}}^2\; dV_{\tilde{g}}\leq C_1(r),
\end{equation}
for some constant $C_1(r)=C_1(K(r),v(r),\bar{\delta}(r))<\infty$. Together with \eqref{eq_cov} this implies
\begin{equation}
 \int_{B_{r/2}(p)}\abs{\widetilde{\Rm}}_{\tilde{g}}^2\; dV_{\tilde{g}}\leq N(r)C_1(r).
\end{equation}
Since $r$ was arbitrary, we can rewrite this as 
\begin{equation}\label{l2tilde}
 \int_{B_{r}(p)}\abs{\widetilde{\Rm}}_{\tilde{g}}^2\; dV_{\tilde{g}}\leq C_2(r),
\end{equation}
where $C_2(r)=N(2r)C_1(2r)$.\\

The final step is to use the $L^2$ bounds \eqref{l2tilde} for $\tilde{g}$ to derive the $L^2$ bounds \eqref{eq_toshow} for the original metric $g$. 
To this end, we first recall the formula for the conformal transformation of the Riemannn tensor, see e.g. \cite{Besse},
\begin{align}
 \Rm=e^{f}\widetilde{\Rm}-g\wedge(\tfrac12\Hess f +\tfrac14\nabla f\otimes \nabla f-\tfrac{1}{8}\abs{\D f}^2)g.
\end{align}
Note that by \eqref{equivmet} it does not matter, up to a factor depending on $r$, whether we compute the norms and volumes with respect to $g$ or $\tilde{g}$. Similarly, by \eqref{lemma1}
the factor $e^f$ can be estimated by a constant depending only on $r$. Thus, the estimate \eqref{l2tilde} implies the estimate \eqref{eq_toshow}, provided that we can estimate the
$L^2$ norms of the terms with the Hessian and the gradient of $f$. By \eqref{2.ballseqn}, \eqref{uppersum} and \eqref{posscal} we have
\begin{equation}
 \int_{B_r(p)}\abs{\D f}^2 dV\leq \omega(r+20)^6.
\end{equation}
By the soliton equation \eqref{1.solitoneq} we can replace $\Hess f$  by $\tfrac{1}{2}g-\Rc$. The $L^2$-norm of the metric is a lower order term and can be easily estimated thanks to \eqref{2.ballseqn}.
For the leading order term, by \cite[Lem. 4.1]{HM11} we have the weighted $L^2$ estimate
\begin{equation}
 \int_M \abs{\Rc}^2 e^{-f} dV\leq C(\ul{\mu}),
\end{equation}
for some constant $C(\ul{\mu})<\infty$. Using again \eqref{lemma1} this implies unweighted $L^2$ bounds for $\Rc$ on $B_r(p)$. Putting everything together, we conclude that
\begin{equation}\label{eq_shown}
 \int_{B_r(p)}\abs{\Rm}_g^2\; dV_g\leq C(r)
\end{equation}
for some universal function $C(r)=C_{\underline{\mu}}(r)<\infty$.\\

Having established the $L^2$-bounds \eqref{eq_shown}, we can now apply our general compactness theorem \cite[Thm 1.1]{HM11}.
This proves Theorem \ref{main_thm}.

\makeatletter
\def\@listi{%
  \itemsep=0pt
  \parsep=1pt
  \topsep=1pt}
\makeatother
{\fontsize{10}{11}\selectfont
}
\vspace{10mm}

Robert Haslhofer\\
{\sc Courant Institute of Mathematical Sciences, New York University, 251 Mercer Street, New York, NY 10012, USA}\\

Reto M\"uller\\
{\sc School of Mathematical Sciences, Queen Mary University of London, Mile End Road, London E1 4NS, UK}\\

email: robert.haslhofer@cims.nyu.edu, r.mueller@qmul.qc.uk

\begin{thebibliography}{99}

\bibitem{Besse}
A.~Besse.
\newblock Einstein manifolds.
\newblock {\em Ergebnisse der Mathematik und ihrer Grenzgebiete (3)}, 10. Springer-Verlag, Berlin, 1987. 

\bibitem{Cao09}
H.-D.~Cao.
\newblock Recent progress on {R}icci solitons.
\newblock In {\em Recent Advances in Geometric Analysis}, volume~11 of {\em
  Advanced Lectures in Mathematics (ALM)}. International Press, 2009.

\bibitem{CS07}
H.-D. Cao and N.~Sesum.
\newblock A compactness result for {K}\"ahler {R}icci solitons.
\newblock {\em Adv. Math.}, 211(2):794--818, 2007.

\bibitem{CZ09}
H.-D. Cao and D.~Zhou.
\newblock On complete gradient shrinking {R}icci solitons.
\newblock {\em J. Differential Geom.}, 85(2):175--185, 2010. 

\bibitem{CN14}
J.~Cheeger and A.~Naber.
\newblock Regularity of {E}instein manifolds and the codimension 4 conjecture.
\newblock {\em arXiv:1406.6534}, 2014.

\bibitem{CW12}
X.~Chen and B.Wang.
\newblock Space of {R}icci flows {I}.
\newblock {\em Comm. Pure Appl. Math.}, 65(10):1399--1457, 2012. 

\bibitem{Ham95s}
R.~Hamilton.
\newblock The formation of singularities in the Ricci flow.
\newblock {\em Surveys in differential geometry}, Vol. II, 7--136, Int. Press, Cambridge, MA, 1995. 

\bibitem{HM11}
R.~Haslhofer and R.~M\"uller.
\newblock A compactness theorem for complete Ricci shrinkers.
\newblock {\em  Geom. Funct. Anal.}, 21(5):1091--1116, 2011.

\bibitem{Per02}
G.~Perelman.
\newblock The entropy formula for the {R}icci flow and its geometric
  applications.
\newblock {\em ArXiv:math/0211159v1}, 2002.

\bibitem{TZ12}
G.~Tian and Z.~Zhang.
\newblock Degeneration of K\"ahler-Ricci solitons.
\newblock {\em Int. Math. Res. Not.} 2012(5):957--985, 2012. 

\bibitem{Web08}
B.~Weber.
\newblock Convergence of compact {R}icci solitons.
\newblock {\em Int. Math. Res. Not.}, 2011(1):96--118, 2011.

\bibitem{Zha06}
X.~Zhang.
\newblock Compactness theorems for gradient {R}icci solitons.
\newblock {\em J. Geom. Phys.}, 56(12):2481--2499, 2006.

\bibitem{Zha09}
Z.H.~Zhang.
\newblock On the completeness of gradient {R}icci solitons.
\newblock {\em Proc. Amer. Math. Soc.}, 137(8):2755--2759, 2009.

\bibitem{Z10}
Z.~Zhang.
\newblock Degeneration of shrinking Ricci solitons. 
\newblock {\em Int. Math. Res. Not.} 2010(21):4137--4158, 2010. 

\end{thebibliography}
\end{document}